\documentclass[12pt]{article}%
\usepackage{amsfonts}
\usepackage{sw20bams}%
\usepackage{amsmath}%
\usepackage{amssymb}%
\usepackage{graphicx}
\providecommand{\U}[1]{\protect\rule{.1in}{.1in}}
\begin{document}

\title{How Far Might We Walk at Random?}
\author{Steven Finch}
\date{February 13, 2018}
\maketitle

\begin{abstract}
This elementary treatment first summarizes extreme values of a Bernoulli
random walk on the one-dimensional integer lattice over a finite discrete time
interval. \ Both the symmetric (unbiased) and asymmetric (biased) cases are
discussed. \ Asymptotic results are given as the time interval length
approaches infinity. \ Focus then shifts to such walks reflected at the origin
-- in both strong and weak senses -- and related unsolved problems are
meticulously examined.

\end{abstract}

\footnotetext{Copyright \copyright \ 2018 by Steven R. Finch. All rights
reserved.}Let $X_{0}=0$ and $X_{1}$, $X_{2}$, \ldots, $X_{n}$ be a sequence of
independent random variables satisfying%
\[%
\begin{array}
[c]{ccccccc}%
\mathbb{P}\left\{  X_{i}=1\right\}  =p, &  & \mathbb{P}\left\{  X_{i}%
=-1\right\}  =q, &  & p+q=1, &  & p\leq q
\end{array}
\]
for each $1\leq i\leq n$. \ Define $S_{0}=X_{0}$ and%
\[
S_{j}=%
{\displaystyle\sum\limits_{i=0}^{j}}
X_{i}=S_{j-1}+X_{j}%
\]
for each $1\leq j\leq n$. \ The simple random walk $S_{0}$, $S_{1}$, $S_{2}$,
\ldots, $S_{n}$ is symmetric if $p=q$ and asymmetric if $p<q$. \ Let%
\[%
\begin{array}
[c]{ccc}%
M_{n}^{+}=\max\limits_{0\leq j\leq n}S_{j}, &  & M_{n}^{-}=-\min\limits_{0\leq
j\leq n}S_{j}%
\end{array}
\]
denote the maximum and (absolute) minimum values of the walk over the time
interval $[0,n]$. \ The quantities $M_{n}^{+}$, $M_{n}^{-}$ are indisputably
the most prominent features of a walk, yet their first/second moments appear
not to be widely known. \ For the symmetric case, the mean of $M_{n}^{-}$ is
the same as the mean of $M_{n}^{+}$. \ Asymptotically, \cite{R1-rw, CM-rw,
Pa-rw, An-rw}%
\[
\mathbb{E}\left(  M_{n}^{+}\right)  \sim\sqrt{\dfrac{2n}{\pi}}-\dfrac{1}{2}%
\]
as $n\rightarrow\infty$. \ The mean square of $M_{n}^{-}$ is the same as the
mean square of $M_{n}^{+}$:%
\[
\mathbb{E}\left(  \left(  M_{n}^{+}\right)  ^{2}\right)  \sim n-\sqrt
{\dfrac{2n}{\pi}}+\dfrac{1}{2}%
\]
which implies that the variance is
\[
\mathbb{V}\left(  M_{n}^{+}\right)  \sim\left(  1-\frac{2}{\pi}\right)  n.
\]
The cross-moment between extremities is \cite{GK-rw, RS-rw, RZ-rw, FP-rw}%
\[
\mathbb{E}\left(  M_{n}^{+}\,M_{n}^{-}\right)  \sim\left(  2\ln(2)-1\right)
n.
\]
An additional term in the asymptotic expansion is not known, perhaps because
calculating $\mathbb{E}\left(  M_{n}^{+}\,M_{n}^{-}\right)  $ seems to
necessitate continuous-time methods \cite{Fr-rw, Cs-rw, KB-rw, CR-rw, Ma-rw}.

For the asymmetric case, $M_{n}^{+}$ follows a geometric distribution as
$n\rightarrow\infty$ (versus a half-normal distribution when $p=1/2$). \ In
particular,%
\[%
\begin{array}
[c]{ccc}%
\mathbb{E}\left(  M_{n}^{+}\right)  \sim\dfrac{p}{1-2p}, &  & \mathbb{E}%
\left(  \left(  M_{n}^{+}\right)  ^{2}\right)  \sim\dfrac{p}{\left(
1-2p\right)  ^{2}}%
\end{array}
\]
are bounded, which is unsurprising since $p<1/2$ and thus the walk has a
downward trend. \ In contrast,%
\[
\mathbb{E}\left(  M_{n}^{-}\right)  \sim(1-2p)n+\frac{p}{1-2p},
\]%
\[
\mathbb{E}\left(  \left(  M_{n}^{-}\right)  ^{2}\right)  \sim(1-2p)^{2}%
n^{2}+2(3-2p)p\,n-\frac{(3-4p)p}{\left(  1-2p\right)  ^{2}},
\]%
\[
\mathbb{E}\left(  M_{n}^{+}\,M_{n}^{-}\right)  \sim p\,n-\frac{(2-3p)p}%
{\left(  1-2p\right)  ^{2}}
\]
are unbounded. \ These results have implications for the range $M_{n}%
^{+}+M_{n}^{-}$ of a random walk \cite{Va-rw, CH-rw, Mo-rw}. \ More will be
said about $M_{n}^{-}$ moments in Section 3.

\section{Strong Reflection at Origin}

Replace the recursive definition $S_{j}=S_{j-1}+X_{j}$ in the preceding
section by%
\[
S_{j}=\left\vert S_{j-1}+X_{j}\right\vert
\]
for all $1\leq j\leq n$. \ It follows that $S_{1}=1$ because $S_{0}=0$ and the
sign of $X_{1}$ becomes immaterial. \ More generally, $S_{k}=1$ whenever
$S_{k-1}=0$. \ Such a \textquotedblleft strong\textquotedblright\ decisive
response to zero crossings is adjusted in the next section. \ Let%
\[
M_{n}=\max\limits_{0\leq j\leq n}S_{j}.
\]
For the symmetric case, Percus \&\ Percus \cite{PP-rw} recently demonstrated
that%
\[%
\begin{array}
[c]{ccc}%
\mathbb{E}\left(  M_{n}\right)  \sim\sqrt{\dfrac{\pi\,n}{2}}, &  &
\mathbb{E}\left(  M_{n}^{2}\right)  \sim2G\,n
\end{array}
\]
as $n\rightarrow\infty$, where%
\[
G=%
{\displaystyle\sum\limits_{k=0}^{\infty}}
\frac{(-1)^{k}}{(2k+1)^{2}}=%
{\displaystyle\int\limits_{0}^{\infty}}
\frac{t}{\cosh(t)}dt=0.9159655941....
\]
is Catalan's constant \cite{F1-rw}. \ These are remarkable formulas! \ No
additional terms in either asymptotic expansion, however, are known.

For the asymmetric case, three procedures are available for generating exact
$M_{n}$ probabilities -- see Section 4 -- and there exists a connection with
enumerating Dyck prefixes of length $n$.

Clearly $M_{n}$ is the largest of independent identically distributed cycle
maxima (cycles occur between every zero crossing of the walk). \ The number of
such cycles cannot exceed $n$. \ Each such maxima possesses a geometric
distribution. \ We conclude that \cite{SR-rw, KhP-rw, Ei-rw}%
\[%
\begin{array}
[c]{ccccc}%
\mathbb{E}\left(  M_{n}\right)  =O\left(  \ln(n)\right)  , &  & \mathbb{E}%
\left(  M_{n}^{2}\right)  =O\left(  \ln(n)^{2}\right)  , &  & \mathbb{V}%
\left(  M_{n}\right)  =O\left(  1\right)
\end{array}
\]
as $n\rightarrow\infty$. \ A\ more precise estimate might involve the actual
(random) cycle count instead of $n$, which brings \cite{SD-rw, Ha-rw} to mind.
\ An alternative argument for the $\mathbb{E}\left(  M_{n}\right)  $ result is
given in Section 5.

\section{Weak Reflection at Origin}

Replace the recursive definition $S_{j}=S_{j-1}+X_{j}$ in the opening section
by%
\[
S_{j}=\max\left\{  S_{j-1}+X_{j},0\right\}
\]
for all $1\leq j\leq n$. \ It follows that $S_{1}=1$ with probability $p$ and
$S_{1}=0$ with probability $q$. \ The same odds apply for $S_{k}$ whenever
$S_{k-1}=0$. \ Such a \textquotedblleft weak\textquotedblright\ anemic
response to zero crossings makes our analysis somewhat more complicated
(relative to our earlier sense of reflection). \ With $M_{n}$ defined exactly
as before, it can be shown that asymptotics
\[%
\begin{array}
[c]{ccc}%
\mathbb{E}\left(  M_{n}\right)  \sim\sqrt{\dfrac{\pi\,n}{2}}, &  &
\mathbb{E}\left(  M_{n}^{2}\right)  \sim2G\,n
\end{array}
\]
likewise hold -- details appear in Section 8 -- and additional terms are again
unknown. \ We conjecture that the discrepancy between strong and weak
reflections should become apparent via higher-order expansions of such moments.

For the asymmetric case, three procedures are available for generating exact
$M_{n}$ probabilities -- see Section 6 -- unfortunately $\mathbb{E}\left(
M_{n}\right)  =O\left(  \ln(n)\right)  $ is again the best known estimate of
the average reflected maximum. \ Elaborate but necessary algebraic details are
given in Section 7.

\section{Recurrence for $(M_{n}^{+},M_{n}^{-})$ Probabilities}

Our methods are non-rigorous and experimentally-based, using the computer
algebra system \textsc{Mathematica}. \ Our starting point is a recurrence for
the joint probability mass function of the maximum $M_{n}^{+}$ and
(absolute)\ minimum $M_{n}^{-}$ of an asymmetric simple random walk
\cite{CxM-rw, RcS-rw, GN-rw}. \ Let $n\geq j$ be positive integers. \ Define
$C[n,j,p,q]$ to be%
\[
\frac{j}{n}\dbinom{n}{h}\dbinom{n-h}{h+j}(p\,q)^{h}
\]
if $n-j$ is even and $0$ otherwise, where $h=(n-j)/2$. \ Given integers $0\leq
a\leq n$ and $0\leq b\leq n$, define $f(n,a,b,p,q)$ to be%
\[
p^{a}%
{\displaystyle\sum\limits_{k=0}^{\ell}}
\left\{  (p\,q)^{(a+b)k}C[n,2(a+b)k+a,p,q]-(p\,q)^{b}%
C[n,2(a+b)k+a+2b,p,q]\right\}
\]
if $n\geq a$ and $0$ otherwise, where $\ell=\left\lfloor (n-a)/\left(
2(a+b)\right)  \right\rfloor $; and $g(n,a,b,p,q)$ to be%
\[
q^{b}%
{\displaystyle\sum\limits_{k=0}^{m}}
\left\{  (p\,q)^{(a+b)k}C[n,2(a+b)k+b,p,q]-(p\,q)^{a}%
C[n,2(a+b)k+b+2a,p,q]\right\}
\]
if $n\geq b$ and $0$ otherwise, where $m=\left\lfloor (n-b)/\left(
2(a+b)\right)  \right\rfloor $. \ Let $\psi(n,a,b,p,q)$ denote%
\[
f(n,a,b,p,q)+g(n,a,b,p,q)
\]
if $n\geq\min\{a,b\}$ and $0$ otherwise. \ The desired recurrence for
\[
\varphi(n,a,b,p,q)=\mathbb{P}\left\{  M_{n}^{+}=a\text{ and }M_{n}%
^{-}=b\right\}
\]
is hence%
\begin{align*}
\varphi(n+1,a,b,p,q)  & =\varphi(n,a,b,p,q)-\psi(n+1,a+1,b+1,p,q)-\psi
(n+1,a,b,p,q)\\
& +\psi(n+1,a+1,b,p,q)+\psi(n+1,a,b+1,p,q)
\end{align*}
with initial/boundary conditions%
\[%
\begin{array}
[c]{ccc}%
\varphi(0,a,b,p,q)=\left\{
\begin{array}
[c]{lll}%
1 &  & \text{if }a=b=0,\\
0 &  & \text{otherwise;}%
\end{array}
\right.  &  & \varphi(n,0,0,p,q)=\left\{
\begin{array}
[c]{lll}%
1 &  & \text{if }n=0,\\
0 &  & \text{otherwise.}%
\end{array}
\right.
\end{array}
\]
For example, the $(n+1)\times(n+1)$ matrix containing entries $\varphi
(n,a,b,p,q)$ for $0\leq a\leq n$ and $0\leq b\leq n$ is%
\[
\left(
\begin{array}
[c]{cc}%
0 & q\\
p & 0
\end{array}
\right)
\]
if $n=1$,
\[
\left(
\begin{array}
[c]{ccc}%
0 & pq & q^{2}\\
pq & 0 & 0\\
p^{2} & 0 & 0
\end{array}
\right)
\]
if $n=2$, and%
\[
\left(
\begin{array}
[c]{cccc}%
0 & pq^{2} & pq^{2} & q^{3}\\
p^{2}q & pq & 0 & 0\\
p^{2}q & 0 & 0 & 0\\
p^{3} & 0 & 0 & 0
\end{array}
\right)
\]
if $n=3$. \ 

Partial verification of the preceding recurrence is provided by formulas
\cite{GS-rw}%
\[
\omega(n,a,p,q)-\omega(n,a+1,p,q)
\]
for the marginal probability mass function of $M_{n}^{+}$ and%
\[
\omega(n,b,q,p)-\omega(n,b+1,q,p)
\]
for the marginal probability mass function of $M_{n}^{-}$, where%
\[
\omega(n,c,x,y)=%
{\displaystyle\sum\limits_{k=c}^{n}}
\left[  1+\left(  \frac{y}{x}\right)  ^{k-c}\right]  \dbinom{n}{\frac{n+k}{2}%
}x^{(n+k)/2}y^{(n-k)/2}.
\]
As before, the binomial coefficient is defined to be $0$ if $n+k$ is odd. Such
consistency checks are helpful; an analogous check of the interaction between
$M_{n}^{+}$, $M_{n}^{-}$ for arbitrary $p$, $q$ is not known.

Series expansions easily indicate that%
\begin{align*}%
{\displaystyle\sum\limits_{a=0}^{n}}
{\displaystyle\sum\limits_{b=0}^{n}}
b\,\varphi(n,a,b,p,1-p)  & \approx n-(2n-1)p+%
{\displaystyle\sum\limits_{i=2}^{n/2}}
2^{i-1}p^{i}\\
& \approx n(1-2p)+%
{\displaystyle\sum\limits_{i=1}^{n/2}}
2^{i-1}p^{i},
\end{align*}%
\begin{align*}%
{\displaystyle\sum\limits_{a=0}^{n}}
{\displaystyle\sum\limits_{b=0}^{n}}
b^{2}\varphi(n,a,b,p,1-p)  & \approx n^{2}-(4n^{2}-6n+3)p+(4n^{2}-4n-8)p^{2}-%
{\displaystyle\sum\limits_{i=3}^{n/2}}
(i+2)2^{i-1}p^{i}\\
& \approx n^{2}(1-4p+4p^{2})+n(6p-4p^{2})-%
{\displaystyle\sum\limits_{i=1}^{n/2}}
(i+2)2^{i-1}p^{i},
\end{align*}%
\begin{align*}%
{\displaystyle\sum\limits_{a=0}^{n}}
{\displaystyle\sum\limits_{b=0}^{n\,}}
a\,b\,\varphi(n,a,b,p,1-p)  & \approx(n-2)p-%
{\displaystyle\sum\limits_{i=2}^{n/3}}
(i+3)2^{i-2}p^{i}\\
& \approx n\,p-%
{\displaystyle\sum\limits_{i=1}^{n/3}}
(i+3)2^{i-2}p^{i}%
\end{align*}
and, after evaluating the sums on the right-hand side as $n\rightarrow\infty$,
the $M_{n}^{-}$ moment asymptotics follow.

With regard to the symmetric case,%
\[%
\begin{array}
[c]{ccc}%
\mathbb{P}\left\{  M_{n}^{+}=k\right\}  =\dfrac{1}{2^{n}}\dbinom
{n}{\left\lfloor \frac{n-k}{2}\right\rfloor }, &  & k=0,1,\ldots,n
\end{array}
\]
thus%
\[%
\begin{array}
[c]{ccc}%
\mathbb{E}\left(  M_{n}^{+}\right)  =\dfrac{\xi_{n}}{2^{n}}, &  &
\mathbb{E}\left(  \left(  M_{n}^{+}\right)  ^{2}\right)  =\dfrac{\eta_{n}%
}{2^{n}}%
\end{array}
\]
where $\xi_{n}$, $\eta_{n}$ are the coefficients of $x^{n}$, $y^{n}$
respectively in the series expansions of \cite{S1-rw}%
\[%
\begin{array}
[c]{ccc}%
\dfrac{-1+2x+\sqrt{1-4x^{2}}}{2\left(  1-2x\right)  ^{2}}, &  & \dfrac
{1+2y-\sqrt{1-4y^{2}}}{2\left(  1-2y\right)  ^{2}}.
\end{array}
\]
The asymptotic results for $\mathbb{E}\left(  M_{n}^{+}\right)  $,
$\mathbb{E}\left(  \left(  M_{n}^{+}\right)  ^{2}\right)  $ are true because%
\[%
\begin{array}
[c]{ccc}%
\xi_{n}\sim\sqrt{\dfrac{2n}{\pi}}2^{n}-2^{n-1}, &  & \eta_{n}=n\,2^{n}-\xi
_{n}\sim n\,2^{n}-\sqrt{\dfrac{2n}{\pi}}2^{n}+2^{n-1}%
\end{array}
\]
as $n\rightarrow\infty$. \ [Note the typographical error for the variance in
\cite{An-rw}, formula (7c).] \ This argument can be made more rigorous, in the
manner of \cite{R1-rw}, without difficulty.

\section{$M_{n}$ Probabilities for Strong Sense}

Henceforth we shall suppress the functional dependence on $p$ and $q$. \ One
procedure for generating exact $M_{n}$ probabilities makes use of an infinite
family of Markov chain matrices:%
\[%
\begin{array}
[c]{ccccc}%
K_{1}=\left(
\begin{array}
[c]{cc}%
0 & 1\\
0 & 1
\end{array}
\right)  , &  & K_{2}=\left(
\begin{array}
[c]{ccc}%
0 & 1 & 0\\
q & 0 & p\\
0 & 0 & 1
\end{array}
\right)  , &  & K_{3}=\left(
\begin{array}
[c]{cccc}%
0 & 1 & 0 & 0\\
q & 0 & p & 0\\
0 & q & 0 & p\\
0 & 0 & 0 & 1
\end{array}
\right)  ,
\end{array}
\]%
\[%
\begin{array}
[c]{ccccc}%
K_{4}=\left(
\begin{array}
[c]{ccccc}%
0 & 1 & 0 & 0 & 0\\
q & 0 & p & 0 & 0\\
0 & q & 0 & p & 0\\
0 & 0 & q & 0 & p\\
0 & 0 & 0 & 0 & 1
\end{array}
\right)  , &  & K_{5}=\left(
\begin{array}
[c]{cccccc}%
0 & 1 & 0 & 0 & 0 & 0\\
q & 0 & p & 0 & 0 & 0\\
0 & q & 0 & p & 0 & 0\\
0 & 0 & q & 0 & p & 0\\
0 & 0 & 0 & q & 0 & p\\
0 & 0 & 0 & 0 & 0 & 1
\end{array}
\right)  , &  & \ldots
\end{array}
\]
Let $\left\{  \varepsilon_{a,1},\varepsilon_{a,2},\ldots,\varepsilon
_{a,a+1}\right\}  $ denote the standard basis of $\mathbb{R}^{a+1}$, then
$\mathbb{P}\left\{  M_{n}=a\right\}  $ is given by \cite{MD-rw, BD-rw, LM-rw}%
\[
\left\{
\begin{array}
[c]{lll}%
1-\varepsilon_{1,1\,}^{^{\prime}}K_{1}^{n}\;\varepsilon_{1,2} &  & \text{if
}a=0,\\
\varepsilon_{a,1\,}^{^{\prime}}K_{a}^{n}\;\varepsilon_{a,a+1}-\varepsilon
_{a+1,1\,}^{^{\prime}}K_{a+1}^{n}\;\varepsilon_{a+1,a+2} &  & \text{if }1\leq
a\leq n
\end{array}
\right.
\]
and $K^{n}$ denotes $n^{\text{th}}$ matrix power. \ Another procedure employs
a recurrence for $\mathbb{P}\left\{  S_{n}=x\text{ and }M_{n}=a\right\}  $:%

\[
F_{1}(x,a)=\delta_{x,1}\delta_{a,1},
\]%
\[
F_{n+1}(0,a)=q\,F_{n}(1,a),
\]%
\[
F_{n+1}(1,a)=F_{n}(0,a)+q\,F_{n}(2,a)\left(  1-\delta_{a,1}\right)  ,
\]%
\[
F_{n+1}(x,a)=\left\{
\begin{array}
[c]{lll}%
p\,F_{n}(x-1,a)+q\,F_{n}(x+1,a) &  & \text{if }1<x<a,\\
p\,F_{n}(a-1,a-1)+p\,F_{n}(a-1,a) &  & \text{if }1<x=a
\end{array}
\right.
\]
and then summing over all $x=0,1,\ldots,a$. \ A\ third procedure is to
determine the coefficient of $\lambda^{n}$ in the series expansion of%
\[
\frac{\theta^{2}-4pq}{1-\lambda}\left[  \frac{2^{a}p^{a-1}\theta^{a}}%
{2^{2a}p^{a-1}q^{a+1}\left(  \theta^{2}-4p^{2}\right)  +\theta^{2a}\left(
\theta^{2}-4q^{2}\right)  }-\frac{2^{a+1}p^{a}\theta^{a+1}}{2^{2a+2}%
p^{a}q^{a+2}\left(  \theta^{2}-4p^{2}\right)  +\theta^{2a+2}\left(  \theta
^{2}-4q^{2}\right)  }\right]
\]
where%
\[
\theta(\lambda)=\dfrac{1-\sqrt{1-4pq\lambda^{2}}}{\lambda}.
\]
The series expansion formula would be less cumbersome upon replacing $\theta$
by $2\theta$, but we wished to retain notational consistency with
\cite{PP-rw}. \ 

Consider an alphabet $\{\nearrow,\searrow\}$. \ A \textbf{Dyck prefix} of
length $n$ is a binary word of a total of $n$ $\nearrow$s and $\searrow$s in
which no initial segment contains more $\searrow$s than $\nearrow$s. \ In the
positive quadrant of the discrete $xy$-plane, starting at the origin,
$\nearrow$ corresponds to the vector $(1,1)$ and $\searrow$ corresponds to the
vector $(1,-1)$. \ Clearly the height of the word (above the $x$-axis) at any
point is nonnegative. \ Counting such words is discussed in \cite{Ke-rw} and
summarized in \cite{S2-rw}. \ The asymptotics given in \cite{Ke-rw} refer to
counts unweighted by probabilities $p$ and $q$; thus they do not help us in
our study.

\section{Alternative $O(\ln(n))$ Argument}

Consider, as in the beginning, an asymmetric simple random walk ($p<q$) that
starts at $0$. \ As in \cite{SD-rw, Ha-rw}, define $T$ to be the time of first
return to $0$, and $M_{T}$ to be the maximum distance from $0$ attained by the
walk up to time $T$ (equivalently, the largest of $M_{T}^{+}$ and $M_{T}^{-}%
$). Conditional on $T<\infty$, we have \cite{Fe-rw}%
\[
\mathbb{P}\left\{  M_{T}<k\right\}  =\frac{(q/p)^{k}-(q/p)}{(q/p)^{k}-1}%
=\frac{(1/x)^{k}-(1/x)}{(1/x)^{k}-1}=\frac{1-x^{k-1}}{1-x^{k}}
\]
where $0<x=p/q<1$, and%
\begin{align*}
\mathbb{P}\left\{  M_{T}=k\right\}   & =\frac{1-x^{k}}{1-x^{k+1}}%
-\frac{1-x^{k-1}}{1-x^{k}}\\
& =\left(  \frac{x^{k-1}}{1-x^{k}}-\frac{x^{k}}{1-x^{k+1}}\right)  +\left(
\frac{1}{1-x^{k+1}}-\frac{1}{1-x^{k}}\right) \\
& =\left(  \frac{x^{k-1}}{1-x^{k}}-\frac{x^{k}}{1-x^{k+1}}\right)  +\left(
\frac{x^{k+1}}{1-x^{k+1}}-\frac{x^{k}}{1-x^{k}}\right) \\
& =\left(  \frac{1}{x}-1\right)  \left(  \frac{x^{k}}{1-x^{k}}-\frac{x^{k+1}%
}{1-x^{k+1}}\right)
\end{align*}
because%
\[
\left(  1-x^{k}\right)  -\left(  1-x^{k+1}\right)  =x^{k+1}\left(
1-x^{k}\right)  -x^{k}\left(  1-x^{k+1}\right)  .
\]
Note that%
\begin{align*}%
{\displaystyle\sum\limits_{k=1}^{\infty}}
k\left(  \frac{x^{k}}{1-x^{k}}-\frac{x^{k+1}}{1-x^{k+1}}\right)   & =%
{\displaystyle\sum\limits_{k=1}^{\infty}}
\frac{k\,x^{k}}{1-x^{k}}-%
{\displaystyle\sum\limits_{k=1}^{\infty}}
\frac{k\,x^{k+1}}{1-x^{k+1}}\\
& =%
{\displaystyle\sum\limits_{k=1}^{\infty}}
\frac{k\,x^{k}}{1-x^{k}}-%
{\displaystyle\sum\limits_{k=2}^{\infty}}
\frac{(k-1)x^{k}}{1-x^{k}}\\
& =%
{\displaystyle\sum\limits_{k=1}^{\infty}}
\frac{k\,x^{k}}{1-x^{k}}-%
{\displaystyle\sum\limits_{k=2}^{\infty}}
\frac{k\,x^{k}}{1-x^{k}}+%
{\displaystyle\sum\limits_{k=2}^{\infty}}
\frac{x^{k}}{1-x^{k}}\\
& =\frac{x}{1-x}+%
{\displaystyle\sum\limits_{k=2}^{\infty}}
\frac{x^{k}}{1-x^{k}}=%
{\displaystyle\sum\limits_{k=1}^{\infty}}
\frac{x^{k}}{1-x^{k}}%
\end{align*}
and therefore%
\[
\mathbb{E}\left(  M_{T}\right)  =\frac{1-x}{x}%
{\displaystyle\sum\limits_{k=1}^{\infty}}
\frac{x^{k}}{1-x^{k}}.
\]
In particular, if $p=1/3$, then $x=1/2$ and \cite{F2-rw, Kn-rw}%
\[
\mathbb{E}\left(  M_{T}\right)  =%
{\displaystyle\sum\limits_{k=1}^{\infty}}
\frac{1}{2^{k}-1}=1.6066951524....
\]
More generally, if $M_{T}^{(n)}$ denotes the largest of $n$ independent copies
of $M_{T}$, we have%
\[
\mathbb{P}\left\{  M_{T}^{(n)}=k\right\}  =\left(  \frac{1-x^{k}}{1-x^{k+1}%
}\right)  ^{n}-\left(  \frac{1-x^{k-1}}{1-x^{k}}\right)  ^{n},
\]%
\[
\mathbb{E}\left(  M_{T}^{(n)}\right)  =%
{\displaystyle\sum\limits_{k=1}^{\infty}}
k\left[  \left(  \frac{1-x^{k}}{1-x^{k+1}}\right)  ^{n}-\left(  \frac
{1-x^{k-1}}{1-x^{k}}\right)  ^{n}\right]  .
\]
For simplicity's sake, set $p=1/3$ (effective for the remainder of this
section), giving%
\begin{align*}
\mathbb{E}\left(  M_{T}^{(n)}\right)   & =%
{\displaystyle\sum\limits_{k=1}^{\infty}}
k\left[  \left(  \frac{2^{k+1}-2}{2^{k+1}-1}\right)  ^{n}-\left(  \frac
{2^{k}-2}{2^{k}-1}\right)  ^{n}\right] \\
& =%
{\displaystyle\sum\limits_{k=1}^{\infty}}
k\left[  \left(  1-\frac{1}{2^{k+1}-1}\right)  ^{n}-\left(  1-\frac{1}%
{2^{k}-1}\right)  ^{n}\right] \\
& =%
{\displaystyle\sum\limits_{k=1}^{\infty}}
k\left\{  \left[  1-\left(  1-\frac{1}{2^{k}-1}\right)  ^{n}\right]  -\left[
1-\left(  1-\frac{1}{2^{k+1}-1}\right)  ^{n}\right]  \right\}  .
\end{align*}
Summation by parts:%
\[%
{\displaystyle\sum\limits_{k=1}^{L}}
u_{k}(v_{k}-v_{k-1})=u_{L}v_{L}-u_{0}v_{0}-%
{\displaystyle\sum\limits_{k=1}^{L}}
(u_{k}-u_{k-1})v_{k-1}
\]
is helpful here with%
\[%
\begin{array}
[c]{ccc}%
u_{k}=k, &  & v_{k}=1-\left(  1-\dfrac{1}{2^{k+1}-1}\right)  ^{n}.
\end{array}
\]
Note that, as $L\rightarrow\infty$,%
\[
u_{L}v_{L}\sim\frac{L\,n}{2^{L+1}-1}\rightarrow0
\]
for any fixed $n$. \ It follows that%
\begin{align*}
\mathbb{E}\left(  M_{T}^{(n)}\right)   & =%
{\displaystyle\sum\limits_{k=1}^{\infty}}
\left[  1-\left(  1-\frac{1}{2^{k}-1}\right)  ^{n}\right] \\
& \sim%
{\displaystyle\sum\limits_{j=0}^{\infty}}
\left[  1-\left(  1-\frac{1}{2^{j}}\right)  ^{n}\right] \\
& \sim\frac{\ln(n)}{\ln(2)}+\left[  \frac{\gamma}{\ln(2)}+\frac{1}{2}+R\left(
\frac{\ln(n)}{\ln(2)}\right)  \right]
\end{align*}
as $n\rightarrow\infty$, where $\gamma$ is Euler's constant \cite{F3-rw} and
$R(z)$ is a periodic function of small amplitude with period $1$ and mean
value $0$. \ This is a famous example, due to Knuth \cite{Kn-rw, SF-rw, VF-rw,
KiP-rw, RJ-rw}, governing the expected cost (number of bit inspectations) per
successful search in a random trie.

Let us address one more calculation: the mean square of $M_{T}$. \ Note that%
\begin{align*}%
{\displaystyle\sum\limits_{k=1}^{\infty}}
k^{2}\left(  \frac{x^{k}}{1-x^{k}}-\frac{x^{k+1}}{1-x^{k+1}}\right)   & =%
{\displaystyle\sum\limits_{k=1}^{\infty}}
\frac{k^{2}\,x^{k}}{1-x^{k}}-%
{\displaystyle\sum\limits_{k=2}^{\infty}}
\frac{(k-1)^{2}\,x^{k}}{1-x^{k}}\\
& =%
{\displaystyle\sum\limits_{k=1}^{\infty}}
\frac{k^{2}\,x^{k}}{1-x^{k}}-%
{\displaystyle\sum\limits_{k=2}^{\infty}}
\frac{k^{2}\,x^{k}}{1-x^{k}}+%
{\displaystyle\sum\limits_{k=2}^{\infty}}
\frac{2k\,x^{k}}{1-x^{k}}-%
{\displaystyle\sum\limits_{k=2}^{\infty}}
\frac{x^{k}}{1-x^{k}}\\
& =\frac{x}{1-x}+\left(
{\displaystyle\sum\limits_{k=1}^{\infty}}
\frac{2k\,x^{k}}{1-x^{k}}-\frac{2x}{1-x}\right)  +\left(  \frac{x}{1-x}-%
{\displaystyle\sum\limits_{k=1}^{\infty}}
\frac{x^{k}}{1-x^{k}}\right) \\
& =2%
{\displaystyle\sum\limits_{k=1}^{\infty}}
\frac{k\,x^{k}}{1-x^{k}}-%
{\displaystyle\sum\limits_{k=1}^{\infty}}
\frac{x^{k}}{1-x^{k}}%
\end{align*}
and therefore%
\[
\mathbb{E}\left(  M_{T}^{2}\right)  =\frac{1-x}{x}\left(  2%
{\displaystyle\sum\limits_{k=1}^{\infty}}
\frac{k\,x^{k}}{1-x^{k}}-%
{\displaystyle\sum\limits_{k=1}^{\infty}}
\frac{x^{k}}{1-x^{k}}\right)  .
\]
In particular, if $p=1/3$, then%
\begin{align*}
\mathbb{E}\left(  M_{T}^{2}\right)   & =2%
{\displaystyle\sum\limits_{k=1}^{\infty}}
\frac{k}{2^{k}-1}-%
{\displaystyle\sum\limits_{k=1}^{\infty}}
\frac{1}{2^{k}-1}\\
& =2(2.7440338887...)-(1.6066951524...)\\
& =3.8813726251....
\end{align*}
Similar work will lead to the asymptotic mean square of $M_{T}^{(n)}$.

\section{$M_{n}$ Probabilities for Weak Sense}

One way for generating exact $M_{n}$ probabilities makes use of Markov chain
matrices:%
\[%
\begin{array}
[c]{ccccc}%
K_{1}=\left(
\begin{array}
[c]{cc}%
q & p\\
0 & 1
\end{array}
\right)  , &  & K_{2}=\left(
\begin{array}
[c]{ccc}%
q & p & 0\\
q & 0 & p\\
0 & 0 & 1
\end{array}
\right)  , &  & K_{3}=\left(
\begin{array}
[c]{cccc}%
q & p & 0 & 0\\
q & 0 & p & 0\\
0 & q & 0 & p\\
0 & 0 & 0 & 1
\end{array}
\right)  ,
\end{array}
\]%
\[%
\begin{array}
[c]{ccccc}%
K_{4}=\left(
\begin{array}
[c]{ccccc}%
q & p & 0 & 0 & 0\\
q & 0 & p & 0 & 0\\
0 & q & 0 & p & 0\\
0 & 0 & q & 0 & p\\
0 & 0 & 0 & 0 & 1
\end{array}
\right)  , &  & K_{5}=\left(
\begin{array}
[c]{cccccc}%
q & p & 0 & 0 & 0 & 0\\
q & 0 & p & 0 & 0 & 0\\
0 & q & 0 & p & 0 & 0\\
0 & 0 & q & 0 & p & 0\\
0 & 0 & 0 & q & 0 & p\\
0 & 0 & 0 & 0 & 0 & 1
\end{array}
\right)  , &  & \ldots
\end{array}
\]
Only two elements (in the northwest corner, first row) are different from
before. \ Evidently $\mathbb{P}\left\{  M_{n}=a\right\}  $ is given by the
identical formula as previously \cite{MD-rw, BD-rw, LM-rw}. \ Another way
employs a recurrence for $\mathbb{P}\left\{  S_{n}=x\text{ and }%
M_{n}=a\right\}  $:%

\[
G_{1}(x,a)=p\,\delta_{x,1}\delta_{a,1}+q\,\delta_{x,0}\delta_{a,0},
\]%
\[
G_{n+1}(0,a)=q\,G_{n}(0,a)+q\,G_{n}(1,a),
\]%
\[
G_{n+1}(1,a)=p\,G_{n}(0,0)\delta_{a,1}+p\,G_{n}(0,a)+q\,G_{n}(2,a)\left(
1-\delta_{a,1}\right)  ,
\]%
\[
G_{n+1}(x,a)=\left\{
\begin{array}
[c]{lll}%
p\,G_{n}(x-1,a)+q\,G_{n}(x+1,a) &  & \text{if }1<x<a,\\
p\,G_{n}(a-1,a-1)+p\,G_{n}(a-1,a) &  & \text{if }1<x=a
\end{array}
\right.
\]
and then summing over all $x=0,1,\ldots,a$. \ A\ third way is to determine the
coefficient of $\lambda^{n}$ in the series expansion of%
\[
\frac{\theta^{2}-4pq}{1-\lambda}\left[  \frac{2^{a}p^{a}\theta^{a}}%
{2^{2a+1}p^{a}q^{a+1}\left(  \theta-2p\right)  +\theta^{2a+1}\left(
\theta-2q\right)  }-\frac{2^{a+1}p^{a+1}\theta^{a+1}}{2^{2a+3}p^{a+1}%
q^{a+2}\left(  \theta-2p\right)  +\theta^{2a+3}\left(  \theta-2q\right)
}\right]
\]
where $\theta(\lambda)$ is the same as before.

Consider an alphabet $\{\nearrow,\rightarrow,\searrow\}$. \ Think of
$\nearrow=(1,1)$ and $\searrow=(1,-1)$ geometrically as before; let
$\rightarrow$ correspond to the horizontal step $(1,0)$. \ A \textbf{dispersed
Dyck prefix} of length $n$ is a ternary word of a total of $n$ $\nearrow$s,
$\rightarrow$s and $\searrow$s in which no initial segment contains more
$\searrow$s than $\nearrow$s; further, no $\rightarrow$ occurs at positive
height. \ Counting such words is summarized in \cite{S2-rw}. \ The associated
(unweighted) asymptotics remain open.

\section{Algebra}

\subsection{Strong Scenario}

Combining elements of the recurrence in Section 4, we have%
\begin{align*}
F_{n+1}(x,a)  & =p\,F_{n}(x-1,a)+q\,\delta_{x,1}F_{n}(x-1,a)-p\,\delta
_{x,a+1}F_{n}(x-1,a)+\\
& q\,F_{n}(x+1,a)+p\,\delta_{x,a}F_{n}(a-1,a-1)\left(  1-\delta_{a,1}\right)
\\
& =\left(  p+q\,\delta_{x,1}-p\,\delta_{x,a+1}\right)  F_{n}(x-1,a)+q\,F_{n}%
(x+1,a)+\\
& p\,\delta_{x,a}F_{n}(a-1,a-1)\left(  1-\delta_{a,1}\right)
\end{align*}
for all $n\geq1$, $0\leq x\leq a$ and $a\geq1$. \ Our goal is to solve for
$F_{n}(x,a)$. \ Introduce the generating function%
\[
F(\lambda,x,a)=%
{\displaystyle\sum\limits_{n=1}^{\infty}}
\lambda^{n}F_{n}(x,a)=\lambda\,F_{1}(x,a)+%
{\displaystyle\sum\limits_{n=1}^{\infty}}
\lambda^{n+1}F_{n+1}(x,a)
\]
from which%
\begin{align*}
F(\lambda,x,a)  & =\lambda\,F_{1}(x,a)+\lambda\left(  p+q\,\delta
_{x,1}-p\,\delta_{x,a+1}\right)  F(\lambda,x-1,a)+\lambda\,q\,F(\lambda
,x+1,a)+\\
& \lambda\,p\,\delta_{x,a}\left(  1-\delta_{a,1}\right)  F(\lambda,a-1,a-1).
\end{align*}
follows. \ Introduce the double generating function\
\[
\tilde{F}(\lambda,\mu,a)=%
{\displaystyle\sum\limits_{n=1}^{\infty}}
{\displaystyle\sum\limits_{x=0}^{a}}
\lambda^{n}\mu^{x}F_{n}(x,a)=%
{\displaystyle\sum\limits_{x=0}^{a}}
F(\lambda,x,a)\mu^{x}
\]
and note that%
\[
\lambda%
{\displaystyle\sum\limits_{x=0}^{a}}
F_{1}(x,a)\mu^{x}=\left\{
\begin{array}
[c]{lll}%
\lambda\,\mu &  & \text{if }a=1\\
0 &  & \text{if }a>1
\end{array}
\right.  \;\;\;=\lambda\,\mu\,\delta_{a,1},
\]%
\[
\lambda\,p%
{\displaystyle\sum\limits_{x=0}^{a}}
F(\lambda,x-1,a)\mu^{x}=\lambda\,p\,\mu%
{\displaystyle\sum\limits_{x=0}^{a-1}}
F(\lambda,x,a)\mu^{x}=\lambda\,p\,\mu\left(  \tilde{F}(\lambda,\mu
,a)-F(\lambda,a,a)\mu^{a}\right)  ,
\]%
\[
\lambda\,q%
{\displaystyle\sum\limits_{x=0}^{a}}
\delta_{x,1}F(\lambda,x-1,a)\mu^{x}-\lambda\,p%
{\displaystyle\sum\limits_{x=0}^{a}}
\delta_{x,a+1}F(\lambda,x-1,a)\mu^{x}=\lambda\,q\,\mu\,F(\lambda,0,a)-0,
\]%
\[
\lambda\,q%
{\displaystyle\sum\limits_{x=0}^{a}}
F(\lambda,x+1,a)\mu^{x}=\,\frac{\lambda\,q}{\mu}%
{\displaystyle\sum\limits_{x=1}^{a+1}}
F(\lambda,x,a)\mu^{x}=\frac{\lambda\,q}{\mu}\left(  \tilde{F}(\lambda
,\mu,a)-F(\lambda,0,a)\mu^{0}\right)  ,
\]%
\[
\lambda\,p%
{\displaystyle\sum\limits_{x=0}^{a}}
\delta_{x,a}\left(  1-\delta_{a,1}\right)  F(\lambda,a-1,a-1)\mu^{x}%
=\lambda\,p\left(  1-\delta_{a,1}\right)  F(\lambda,a-1,a-1)\mu^{a}.
\]
We obtain%
\begin{align*}
\tilde{F}  & =\lambda\,\mu\,\delta_{a,1}+\lambda\,p\,\mu\left(  \tilde
{F}-F(\lambda,a,a)\mu^{a}\right)  +\lambda\,q\,\mu\,F(\lambda,0,a)\\
& +\frac{\lambda\,q}{\mu}\left(  \tilde{F}-F(\lambda,0,a)\right)
+\lambda\,p\left(  1-\delta_{a,1}\right)  F(\lambda,a-1,a-1)\mu^{a}%
\end{align*}
that is,%
\begin{align*}
\left(  1-\lambda\,p\,\mu-\frac{\lambda\,q}{\mu}\right)  \tilde{F}  &
=\lambda\,\mu\,\delta_{a,1}+\left(  \lambda\,q\,\mu-\frac{\lambda\,q}{\mu
}\right)  F(\lambda,0,a)-\lambda\,p\,\mu^{a+1}F(\lambda,a,a)+\\
& \lambda\,p\left(  1-\delta_{a,1}\right)  \mu^{a}F(\lambda,a-1,a-1)
\end{align*}
that is,%
\begin{align}
\left(  p\,\mu^{2}-\frac{\mu}{\lambda}+q\right)  \tilde{F}  & =-\mu^{2}%
\delta_{a,1}+q\left(  1-\mu^{2}\right)  F(\lambda,0,a)+p\,\mu^{a+2}%
F(\lambda,a,a)-\nonumber\\
& p\left(  1-\delta_{a,1}\right)  \mu^{a+1}F(\lambda,a-1,a-1)
\end{align}
after multiplying both sides by $-\mu/\lambda$. \ Examine the special case
$a=1$:%
\[
p\left(  \mu-\frac{\theta}{2p}\right)  \left(  \mu-\frac{2q}{\theta}\right)
\tilde{F}=-\mu^{2}+q\left(  1-\mu^{2}\right)  F(\lambda,0,1)+p\,\mu
^{3}F(\lambda,1,1)
\]
where%
\[%
\begin{array}
[c]{ccccc}%
\dfrac{\theta}{2p}=\dfrac{1-\sqrt{1-4pq\lambda^{2}}}{2p\lambda}, &  &
\text{equivalently,} &  & \dfrac{2q}{\theta}=\dfrac{1+\sqrt{1-4pq\lambda^{2}}%
}{2p\lambda}.
\end{array}
\]
For future reference, the sum of the two zeroes is $1/(p\lambda)$, which
implies that
\[
\lambda=\frac{1}{p}\frac{1}{\tfrac{\theta}{2p}+\tfrac{2q}{\theta}}%
=\frac{2\theta}{\theta^{2}+4pq}.
\]
Taking $\mu=\theta/(2p)$ and then $\mu=2q/\theta$, we have%
\[
\left\{
\begin{array}
[c]{l}%
0=-\dfrac{\theta^{2}}{4p^{2}}+q\left(  1-\dfrac{\theta^{2}}{4p^{2}}\right)
F(\lambda,0,1)+\dfrac{\theta^{3}}{8p^{2}}F(\lambda,1,1)\\
0=-\dfrac{4q^{2}}{\theta^{2}}+q\left(  1-\dfrac{4q^{2}}{\theta^{2}}\right)
F(\lambda,0,1)+\dfrac{8pq^{3}}{\theta^{3}}F(\lambda,1,1)
\end{array}
\right.
\]
and, on eliminating $F(\lambda,0,1)$,
\[
F(\lambda,1,1)=\frac{2^{1}p^{0}\left(  \theta^{2}-4pq\right)  \left(
\theta^{2}+4pq\right)  \theta^{1}}{2^{4}p^{1}q^{3}\left(  \theta^{2}%
-4p^{2}\right)  +\theta^{4}\left(  \theta^{2}-4q^{2}\right)  }.
\]
Now examine the general case $a>1$:%
\[
\left\{
\begin{array}
[c]{l}%
0=q\left(  1-\dfrac{\theta^{2}}{4p^{2}}\right)  F(\lambda,0,a)+\dfrac
{\theta^{a+2}}{2^{a+2}p^{a+1}}F(\lambda,a,a)-\dfrac{\theta^{a+1}}{2^{a+1}%
p^{a}}F(\lambda,a-1,a-1)\\
0=q\left(  1-\dfrac{4q^{2}}{\theta^{2}}\right)  F(\lambda,0,a)+\dfrac
{2^{a+2}pq^{a+2}}{\theta^{a+2}}F(\lambda,a,a)-\dfrac{2^{a+1}pq^{a+1}}%
{\theta^{a+1}}F(\lambda,a-1,a-1)
\end{array}
\right.
\]
and, on eliminating $F(\lambda,0,a)$,%
\begin{align*}
F(\lambda,a,a)  & =2p\frac{2^{2a}p^{a-1}q^{a+1}\left(  \theta^{2}%
-4p^{2}\right)  +\theta^{2a}\left(  \theta^{2}-4q^{2}\right)  }{2^{2a+2}%
p^{a}q^{a+2}\left(  \theta^{2}-4p^{2}\right)  +\theta^{2a+2}\left(  \theta
^{2}-4q^{2}\right)  }\theta\,F(\lambda,a-1,a-1)\\
& =\frac{2^{a}p^{a-1}\left(  \theta^{2}-4pq\right)  \left(  \theta
^{2}+4pq\right)  \theta^{a}}{2^{2a+2}p^{a}q^{a+2}\left(  \theta^{2}%
-4p^{2}\right)  +\theta^{2a+2}\left(  \theta^{2}-4q^{2}\right)  }%
\end{align*}
after iteration. \ Finally, given $a>1$ and taking the limit in formula (1) as
$\mu\rightarrow1$, we have%
\[
-\frac{1-\lambda}{\lambda}\tilde{F}=\left(  1-\frac{1}{\lambda}\right)
\tilde{F}=\left(  p-\frac{1}{\lambda}+q\right)  \tilde{F}=p\,F(\lambda
,a,a)-p\,F(\lambda,a-1,a-1)
\]
therefore%
\begin{align*}
\tilde{F}  & =\frac{p\lambda}{1-\lambda}\left[  F(\lambda,a-1,a-1)-F(\lambda
,a,a)\right] \\
& =\frac{1}{1-\lambda}\frac{2p\theta}{\theta^{2}+4pq}\left[  F(\lambda
,a-1,a-1)-F(\lambda,a,a)\right]
\end{align*}
as was to be shown. \ The case $a=1$ must be treated separately:%
\begin{align*}
\tilde{F}  & =F(\lambda,0,1)+F(\lambda,1,1)=\frac{4p^{2}}{4p^{2}-\theta^{2}%
}\frac{1}{q}\left[  \frac{\theta^{2}}{4p^{2}}-\frac{\theta^{3}}{8p^{2}%
}F(\lambda,1,1)\right]  +F(\lambda,1,1)\\
& =\frac{\lambda\left(  1+q\lambda\right)  }{1-q\lambda^{2}}%
\end{align*}
consistent with the series expansion in Section 4.

\subsection{Weak Scenario}

Combining elements of the recurrence in Section 6, we have%
\begin{align*}
G_{n+1}(x,a)  & =p\,G_{n}(x-1,a)+q\,\delta_{x,0}G_{n}(x,a)-p\,\delta
_{x,a+1}G_{n}(x-1,a)+\\
& q\,G_{n}(x+1,a)+p\,G_{n}(0,0)\delta_{x,1}\delta_{a,1}+p\,\delta_{x,a}%
G_{n}(a-1,a-1)\left(  1-\delta_{a,1}\right) \\
& =\left(  p-p\,\delta_{x,a+1}\right)  G_{n}(x-1,a)+q\,\delta_{x,0}%
G_{n}(x,a)+q\,G_{n}(x+1,a)+\\
& p\,\delta_{x,a}G_{n}(a-1,a-1)
\end{align*}
for all $n\geq1$, $0\leq x\leq a$ and $a\geq1$. \ In particular,
$G_{1}(0,0)=q$, hence $G_{2}(0,0)=q\,G_{1}(0,0)+0=q^{2}$ and $G_{n}%
(0,0)=q\,G_{n-1}(0,0)+0=q^{n}$. \ Our goal is to solve for $G_{n}(x,a)$. \ As
before, we introduce $G(\lambda,x,a)$ and observe that%
\[
G(\lambda,0,0)=%
{\displaystyle\sum\limits_{n=1}^{\infty}}
\lambda^{n}G_{n}(0,0)=%
{\displaystyle\sum\limits_{n=1}^{\infty}}
\lambda^{n}q^{n}=\frac{\lambda\,q}{1-\lambda\,q},
\]%
\begin{align*}
G(\lambda,x,a)  & =\lambda\,G_{1}(x,a)+\lambda\left(  p-p\,\delta
_{x,a+1}\right)  G(\lambda,x-1,a)+\lambda\,q\,\delta_{x,0}G(\lambda,x,a)+\\
& \lambda\,q\,G(\lambda,x+1,a)+\lambda\,p\,\delta_{x,a}G(\lambda,a-1,a-1).
\end{align*}
Introduce $\tilde{G}(\lambda,\mu,a)$ as before and note that%
\[
\lambda%
{\displaystyle\sum\limits_{x=0}^{a}}
G_{1}(x,a)\mu^{x}=\left\{
\begin{array}
[c]{lll}%
\lambda\,q &  & \text{if }a=0\\
\lambda\,p\,\mu &  & \text{if }a=1\\
0 &  & \text{if }a>1
\end{array}
\right.  \;\;\;=\lambda\,q\,\delta_{a,0}+\lambda\,p\,\mu\,\delta_{a,1},
\]%
\[
-p\,\mu^{2}-p\,\mu^{2}G(\lambda,0,0)=-p\,\mu^{2}-\frac{\lambda\,\mu^{2}%
\,p\,q}{1-\lambda\,q}=\frac{-p\,\mu^{2}}{1-\lambda\,q}.
\]
We obtain%
\begin{align*}
\tilde{G}  & =\left(  \lambda\,q\,\delta_{a,0}+\lambda\,p\,\mu\,\delta
_{a,1}\right)  +\lambda\,p\,\mu\left(  \tilde{G}-G(\lambda,a,a)\mu^{a}\right)
+\lambda\,q\,G(\lambda,0,a)\\
& +\frac{\lambda\,q}{\mu}\left(  \tilde{G}-G(\lambda,0,a)\right)
+\lambda\,p\,G(\lambda,a-1,a-1)\mu^{a}%
\end{align*}
that is,%
\begin{align*}
\left(  1-\lambda\,p\,\mu-\frac{\lambda\,q}{\mu}\right)  \tilde{G}  & =\left(
\lambda\,q\,\delta_{a,0}+\lambda\,p\,\mu\,\delta_{a,1}\right)  -\lambda
\,q\left(  \frac{1}{\mu}-1\right)  G(\lambda,0,a)-\\
& \lambda\,p\,\mu^{a+1}G(\lambda,a,a)+\lambda\,p\,\mu^{a}G(\lambda,a-1,a-1)
\end{align*}
that is,%
\begin{align}
\left(  p\,\mu^{2}-\frac{\mu}{\lambda}+q\right)  \tilde{G}  & =-q\,\mu
\,\delta_{a,0}-p\,\mu^{2}\delta_{a,1}+q\left(  1-\mu\right)  G(\lambda,0,a)+\\
& p\,\mu^{a+2}G(\lambda,a,a)-p\,\mu^{a+1}G(\lambda,a-1,a-1)\nonumber
\end{align}
after multiplying both sides by $-\mu/\lambda$. \ Examine the special case
$a=1$:%
\[
p\left(  \mu-\frac{\theta}{2p}\right)  \left(  \mu-\frac{2q}{\theta}\right)
\tilde{G}=-\frac{p\,\mu^{2}}{1-q\,\lambda}+q\left(  1-\mu\right)
G(\lambda,0,1)+p\,\mu^{3}G(\lambda,1,1)
\]
where $\theta$ is exactly as before. \ Taking $\mu=\theta/(2p)$ and then
$\mu=2q/\theta$, we have%
\[
\left\{
\begin{array}
[c]{l}%
0=-\dfrac{\theta^{2}}{4p(1-q\lambda)}+q\left(  1-\dfrac{\theta}{2p}\right)
G(\lambda,0,1)+\dfrac{\theta^{3}}{8p^{2}}G(\lambda,1,1)\\
0=-\dfrac{4pq^{2}}{\theta^{2}(1-q\lambda)}+q\left(  1-\dfrac{2q}{\theta
}\right)  G(\lambda,0,1)+\dfrac{8pq^{3}}{\theta^{3}}G(\lambda,1,1)
\end{array}
\right.
\]
and, on eliminating $G(\lambda,0,1)$,
\[
G(\lambda,1,1)=\frac{2^{1}p^{1}\left(  \theta^{2}-4pq\right)  \left(
\theta^{2}+4pq\right)  \theta^{1}}{2^{5}p^{2}q^{3}\left(  \theta-2p\right)
+\theta^{5}\left(  \theta-2q\right)  }.
\]
Now examine the general case $a>1$:%
\[
\left\{
\begin{array}
[c]{l}%
0=q\left(  1-\dfrac{\theta}{2p}\right)  G(\lambda,0,a)+\dfrac{\theta^{a+2}%
}{2^{a+2}p^{a+1}}G(\lambda,a,a)-\dfrac{\theta^{a+1}}{2^{a+1}p^{a}}%
G(\lambda,a-1,a-1)\\
0=q\left(  1-\dfrac{2q}{\theta}\right)  G(\lambda,0,a)+\dfrac{2^{a+2}pq^{a+2}%
}{\theta^{a+2}}G(\lambda,a,a)-\dfrac{2^{a+1}pq^{a+1}}{\theta^{a+1}}%
G(\lambda,a-1,a-1)
\end{array}
\right.
\]
and, on eliminating $G(\lambda,0,a)$,%
\begin{align*}
G(\lambda,a,a)  & =2p\frac{2^{2a+1}p^{a}q^{a+1}\left(  \theta-2p\right)
+\theta^{2a+1}\left(  \theta-2q\right)  }{2^{2a+3}p^{a+1}q^{a+2}\left(
\theta-2p\right)  +\theta^{2a+3}\left(  \theta-2q\right)  }\theta
\,G(\lambda,a-1,a-1)\\
& =\frac{2^{a}p^{a}\left(  \theta^{2}-4pq\right)  \left(  \theta
^{2}+4pq\right)  \theta^{a}}{2^{2a+3}p^{a+1}q^{a+2}\left(  \theta-2p\right)
+\theta^{2a+3}\left(  \theta-2q\right)  }%
\end{align*}
after iteration. \ It is interesting to compare the expressions for
$F(\lambda,a,a)$ and $G(\lambda,a,a)$. \ Finally, given $a>1$ and taking the
limit in formula (2) as $\mu\rightarrow1$, we have%
\[
-\frac{1-\lambda}{\lambda}\tilde{G}=p\,G(\lambda,a,a)-p\,G(\lambda,a-1,a-1)
\]
therefore%
\[
\tilde{G}=\frac{1}{1-\lambda}\frac{2p\theta}{\theta^{2}+4pq}\left[
G(\lambda,a-1,a-1)-G(\lambda,a,a)\right]
\]
as was to be shown. \ The case $a=1$ again must be treated separately:%
\begin{align*}
\tilde{G}  & =G(\lambda,0,1)+G(\lambda,1,1)=\frac{2p}{2p-\theta}\frac{1}%
{q}\left[  \frac{\theta^{2}}{4p\left(  1-q\lambda\right)  }-\frac{\theta^{3}%
}{8p^{2}}G(\lambda,1,1)\right]  +G(\lambda,1,1)\\
& =\frac{p\lambda}{\left(  1-q\lambda\right)  \left(  1-q\lambda-pq\lambda
^{2}\right)  }%
\end{align*}
consistent with the series expansion in Section 6.

\section{Calculus}

\subsection{Strong Scenario}

Setting $p=q=1/2$, we obtain%
\begin{align*}
\tilde{F}(\lambda,1,a)  & =\frac{2}{1-\lambda}\left(  \frac{\theta^{a}%
}{1+\theta^{2a}}-\frac{\theta^{a+1}}{1+\theta^{2a+2}}\right) \\
& =\frac{2}{1-\lambda}\left(  \frac{1}{\theta^{a}+\theta^{-a}}-\frac{1}%
{\theta^{a+1}+\theta^{-a-1}}\right)
\end{align*}
and thus have the double generating function%
\begin{align*}
\Phi(\lambda,\nu)  & =%
{\displaystyle\sum\limits_{n=1}^{\infty}}
{\displaystyle\sum\limits_{a=1}^{\infty}}
\lambda^{n}\nu^{a}\mathbb{P}\left\{  M_{n}=a\right\}  =%
{\displaystyle\sum\limits_{a=1}^{\infty}}
\nu^{a}\tilde{F}(\lambda,1,a)\\
& =\frac{1}{1-\lambda}%
{\displaystyle\sum\limits_{a=1}^{\infty}}
\nu^{a}\left(  \frac{2}{\theta^{a}+\theta^{-a}}-\frac{2}{\theta^{a+1}%
+\theta^{-a-1}}\right) \\
& =\frac{1}{1-\lambda}\left[  \nu^{0}\frac{2}{\theta^{1}+\theta^{-1}}+%
{\displaystyle\sum\limits_{a=1}^{\infty}}
\left(  \nu^{a}-\nu^{a-1}\right)  \frac{2}{\theta^{a}+\theta^{-a}}\right] \\
& =\frac{\lambda}{1-\lambda}-\frac{1-\nu}{1-\lambda}%
{\displaystyle\sum\limits_{a=1}^{\infty}}
\nu^{a-1}\frac{2}{\theta^{a}+\theta^{-a}}.
\end{align*}
[Note that $\nu^{a-1}$ in the final sum is mistakenly given as $\nu^{a}$ in
\cite{PP-rw}, formula (4.4)]. \ Let us focus on $\mathbb{E}\left(
M_{n}\right)  $ solely:%
\begin{align*}%
{\displaystyle\sum\limits_{n=1}^{\infty}}
\lambda^{n}\mathbb{E}\left(  M_{n}\right)   & =\left.  \frac{\partial\Phi
}{\partial\nu}\right\vert _{\nu=1}\\
& =\left.  \frac{1}{1-\lambda}%
{\displaystyle\sum\limits_{a=1}^{\infty}}
\nu^{a-1}\frac{2}{\theta^{a}+\theta^{-a}}-\frac{1-\nu}{1-\lambda}%
{\displaystyle\sum\limits_{a=1}^{\infty}}
(a-1)\nu^{a-2}\frac{2}{\theta^{a}+\theta^{-a}}\right\vert _{\nu=1}\\
& =\frac{1}{1-\lambda}%
{\displaystyle\sum\limits_{a=1}^{\infty}}
\frac{2}{\theta^{a}+\theta^{-a}}%
\end{align*}
which provides that (in an extended sense) the following sequence is Abel
convergent \cite{PP-rw, PS-rw}:%
\begin{align*}
\lim_{n\rightarrow\infty}^{\;\;\;\;\;\;\;\;\ast}\frac{\mathbb{E}\left(
M_{n}\right)  }{\sqrt{n}}  & =\lim_{\lambda\rightarrow1^{-}}(1-\lambda)^{3/2}%
{\displaystyle\sum\limits_{n=1}^{\infty}}
\frac{1}{(1/2)!}\lambda^{n-1/2}\mathbb{E}\left(  M_{n}\right) \\
& =\lim_{\lambda\rightarrow1^{-}}\left(  \frac{1-\lambda}{\lambda}\right)
^{1/2}%
{\displaystyle\sum\limits_{a=1}^{\infty}}
\frac{2}{\theta^{a}+\theta^{-a}}\frac{1}{(1/2)!}.
\end{align*}
Let $\theta=\exp(-t)$, then%
\[
\frac{2}{\theta^{a}+\theta^{-a}}=\frac{2}{e^{at}+e^{-at}}=\operatorname{sech}%
(at)
\]
and, because $\lambda=2\theta/(\theta^{2}+1)$,
\[
\dfrac{1-\lambda}{\lambda}=\frac{\theta^{2}+1}{2\theta}-1=\frac{\theta
+\theta^{-1}}{2}-1=\cosh(t)-1.
\]
We have%
\begin{align*}
\lim_{n\rightarrow\infty}^{\;\;\;\;\;\;\;\;\ast}\frac{\mathbb{E}\left(
M_{n}\right)  }{\sqrt{n}}  & =\lim_{t\rightarrow0^{+}}\sqrt{\cosh(t)-1}%
\frac{1}{(1/2)!}%
{\displaystyle\sum\limits_{a=1}^{\infty}}
\operatorname{sech}(at)\\
& =\lim_{t\rightarrow0^{+}}\sqrt{\frac{2}{\pi}}t%
{\displaystyle\sum\limits_{a=1}^{\infty}}
\operatorname{sech}(at)
\end{align*}
since $\cosh(t)-1\sim t^{2}/2$ and $(1/2)!=\sqrt{\pi}/2$. \ By a Riemann
sum-based argument,
\[
t%
{\displaystyle\sum\limits_{a=\left\lceil \alpha/t\right\rceil }^{\left\lfloor
\beta/t\right\rfloor }}
\operatorname{sech}(at)\rightarrow%
{\displaystyle\int\limits_{\alpha}^{\beta}}
\operatorname{sech}(b)db
\]
as $t\rightarrow0^{+}$, which in turn gives%
\[
\lim_{n\rightarrow\infty}^{\;\;\;\;\;\;\;\;\ast}\frac{\mathbb{E}\left(
M_{n}\right)  }{\sqrt{n}}=\sqrt{\frac{2}{\pi}}%
{\displaystyle\int\limits_{0}^{\infty}}
\operatorname{sech}(b)db=\sqrt{\frac{\pi}{2}}
\]
as $\alpha\rightarrow0^{+}$ and $\beta\rightarrow\infty$. \ A similar argument
gives the $\mathbb{E}\left(  M_{n}^{2}\right)  $ result in terms of Catalan's constant.

\subsection{Weak Scenario}

The preceding was old (imitating \cite{PP-rw}); the following is apparently
new. \ Setting $p=q=1/2$, we obtain%
\begin{align*}
\tilde{G}(\lambda,1,a)  & =\frac{1+\theta}{1-\lambda}\left[  \frac{\theta^{a}%
}{1+\theta^{2a+1}}-\frac{\theta^{a+1}}{1+\theta^{2a+3}}\right] \\
& =\frac{1}{1-\lambda}\frac{1+\theta}{\sqrt{\theta}}\left[  \frac{1}%
{\theta^{a+1/2}+\theta^{-a-1/2}}-\frac{1}{\theta^{a+3/2}+\theta^{-a-3/2}%
}\right]
\end{align*}
for $a\geq1$ and%
\[
\tilde{G}(\lambda,1,0)=G(\lambda,0,0)=%
{\displaystyle\sum\limits_{n=1}^{\infty}}
\left(  \frac{\lambda}{2}\right)  ^{n}=\frac{\frac{\lambda}{2}}{1-\frac
{\lambda}{2}}.
\]
As a preliminary step,%
\[
\frac{\theta^{1/2}+\theta^{-1/2}}{\theta^{3/2}+\theta^{-3/2}}=\frac{\theta
^{2}+\theta}{\theta^{3}+1}=\frac{\theta}{\theta^{2}-\theta+1}=\frac
{\frac{\theta}{\theta^{2}+1}}{1-\frac{\theta}{\theta^{2}+1}}=\frac
{\frac{\lambda}{2}}{1-\frac{\lambda}{2}}
\]
and%
\[
\frac{\frac{\lambda}{2}}{1-\frac{\lambda}{2}}+\frac{\frac{\lambda}{2}}%
{1-\frac{\lambda}{2}}\frac{1}{1-\lambda}=\frac{\lambda}{2-\lambda}\left(
1+\frac{1}{1-\lambda}\right)  =\frac{\lambda}{2-\lambda}\frac{2-\lambda
}{1-\lambda}=\frac{\lambda}{1-\lambda}.
\]
We thus have%
\begin{align*}
\Psi(\lambda,\nu)  & =%
{\displaystyle\sum\limits_{a=0}^{\infty}}
\nu^{a}\tilde{G}(\lambda,1,a)\\
& =\frac{\lambda}{2-\lambda}+\frac{1}{1-\lambda}\frac{1+\theta}{\sqrt{\theta}}%
{\displaystyle\sum\limits_{a=1}^{\infty}}
\nu^{a}\left(  \frac{1}{\theta^{a+1/2}+\theta^{-a-1/2}}-\frac{1}%
{\theta^{a+3/2}+\theta^{-a-3/2}}\right) \\
& =\frac{\lambda}{2-\lambda}+\frac{1+\theta}{\sqrt{\theta}}\frac{1}{1-\lambda
}\left[  \nu^{0}\frac{1}{\theta^{3/2}+\theta^{-3/2}}+%
{\displaystyle\sum\limits_{a=1}^{\infty}}
\left(  \nu^{a}-\nu^{a-1}\right)  \frac{1}{\theta^{a+1/2}+\theta^{-a-1/2}%
}\right] \\
& =\frac{\lambda}{1-\lambda}-\frac{1+\theta}{\sqrt{\theta}}\frac{1-\nu
}{1-\lambda}%
{\displaystyle\sum\limits_{a=1}^{\infty}}
\nu^{a-1}\frac{1}{\theta^{a+1/2}+\theta^{-a-1/2}}.
\end{align*}
Obtaining $\mathbb{E}\left(  M_{n}\right)  $ is achieved by $\nu
$-differentiation of $\Psi$ at $\nu=1$:
\begin{align*}
& \frac{1+\theta}{\sqrt{\theta}}\left(  \left.  \frac{1}{1-\lambda}%
{\displaystyle\sum\limits_{a=1}^{\infty}}
\nu^{a-1}\frac{1}{\theta^{a+1/2}+\theta^{-a-1/2}}-\frac{1-\nu}{1-\lambda}%
{\displaystyle\sum\limits_{a=1}^{\infty}}
(a-1)\nu^{a-2}\frac{1}{\theta^{a+1/2}+\theta^{-a-1/2}}\right\vert _{\nu
=1}\right) \\
& =\frac{1+\theta}{2\sqrt{\theta}}\frac{1}{1-\lambda}%
{\displaystyle\sum\limits_{a=1}^{\infty}}
\frac{2}{\theta^{a+1/2}+\theta^{-a-1/2}}%
\end{align*}
which provides that the following sequence is Abel convergent:%
\begin{align*}
\lim_{n\rightarrow\infty}^{\;\;\;\;\;\;\;\;\ast}\frac{\mathbb{E}\left(
M_{n}\right)  }{\sqrt{n}}  & =\lim_{\lambda\rightarrow1^{-}}\frac{1+\theta
}{2\sqrt{\theta}}\left(  \frac{1-\lambda}{\lambda}\right)  ^{1/2}%
{\displaystyle\sum\limits_{a=1}^{\infty}}
\frac{2}{\theta^{a+1/2}+\theta^{-a-1/2}}\frac{1}{(1/2)!}\\
& =\lim_{t\rightarrow0^{+}}\frac{1+e^{-t}}{2e^{-t/2}}\sqrt{\cosh(t)-1}\frac
{1}{(1/2)!}%
{\displaystyle\sum\limits_{a=1}^{\infty}}
\operatorname{sech}\left(  \left(  a+\frac{1}{2}\right)  t\right) \\
& =\lim_{t\rightarrow0^{+}}\sqrt{\frac{2}{\pi}}t%
{\displaystyle\sum\limits_{a=1}^{\infty}}
\operatorname{sech}\left(  \left(  a+\frac{1}{2}\right)  t\right)  .
\end{align*}
By a Riemann sum-based argument,
\[
t%
{\displaystyle\sum\limits_{a=\left\lceil \alpha/t\right\rceil }^{\left\lfloor
\beta/t\right\rfloor }}
\operatorname{sech}\left(  \left(  a+\frac{1}{2}\right)  t\right)  \rightarrow%
{\displaystyle\int\limits_{\alpha}^{\beta}}
\operatorname{sech}(b)db
\]
as $t\rightarrow0^{+}$, thus the same limiting constants (as $\alpha
\rightarrow0^{+}$ and $\beta\rightarrow\infty$) for $\mathbb{E}\left(
M_{n}\right)  $ and $\mathbb{E}\left(  M_{n}^{2}\right)  $ apply here as before.

\section{Closing Words}

After the preceding was written, I\ realized that the constants $\sqrt{\pi/2}$
and $2G$ appeared long before Percus \&\ Percus \cite{PP-rw} in a
closely-related context.\ Given a symmetric simple random walk $S_{0}$,
$S_{1}$, $S_{2}$, \ldots, $S_{n}$ as at the begining, the quantity%
\[
M_{n}=\max\{M_{n}^{+},M_{n}^{-}\}=\max\limits_{0\leq j\leq n}\,\left\vert
S_{j}\right\vert
\]
possesses a well-known distribution function due to Erd\"{o}s \&\ Kac
\cite{EK-rw} with the aforementioned moments given explicitly by R\'{e}nyi
\cite{R2-rw}. \ It is not surprising that iterating $S_{j}\leftarrow
S_{j-1}+X_{j} $ and maximizing $\left\vert S_{j}\right\vert $ should be
asymptotically analogous to iterating $S_{j}\leftarrow\left\vert S_{j-1}%
+X_{j}\right\vert $ and maximizing $S_{j}$ (strong reflection). \ Our humble
contribution was to show that iterating $S_{j}\leftarrow\max\left\{
S_{j-1}+X_{j},0\right\}  $ and maximizing $S_{j}$ (weak reflection) likewise
falls in the same realm.

As soon as symmetry is discarded, concrete theory underlying reflected random
walks becomes rather thin. The problem of determining precisely the $O$
constants, as functions of $p<q$, in%
\[%
\begin{array}
[c]{ccc}%
\mathbb{E}\left(  M_{n}\right)  =O\left(  \ln(n)\right)  , &  & \mathbb{E}%
\left(  M_{n}^{2}\right)  =O\left(  \ln(n)^{2}\right)
\end{array}
\]
is apparently open \cite{A1-rw, A2-rw}. \ Doing this in such a way as to
deduce a formula for $\mathbb{V}(M_{n})$, based on series expansion results in
Sections 4 \&\ 6, would seem to be very hard.

If we allow zero to be a third possible value in each of the steps:
\[%
\begin{array}
[c]{ccccccc}%
\mathbb{P}\left\{  X_{i}=1\right\}  =p, &  & \mathbb{P}\left\{  X_{i}%
=-1\right\}  =q, &  & \mathbb{P}\left\{  X_{i}=0\right\}  =r, &  & p+q+r=1
\end{array}
\]
for $1\leq i\leq n$, then there emerges a \textbf{lazy random walk} $S_{0}$,
$S_{1}$, $S_{2}$, \ldots, $S_{n}$ with more unanswered questions. The
recurrence given in Section 3 for joint probabilities in the asymmetric case
can be easily extended \cite{CxM-rw, RcS-rw, GN-rw}. \ Katzenbeisser \& Panny
\cite{KaP-rw} studied the symmetric case ($p=q$) under the extra condition
that $S_{n}=S_{0}$; Prodinger \&\ Panny \cite{P1-rw, P2-rw, P3-rw} studied
likewise but further assumed paths to be nonnegative. \ Their results align
with known formulas for Brownian bridge and Brownian excursion, respectively.
\ We infer that%
\[%
\begin{array}
[c]{ccc}%
\mathbb{E}\left(  M_{n}\right)  \sim\sqrt{\dfrac{\pi\left(  1-r\right)  n}{2}%
}, &  & \mathbb{E}\left(  M_{n}^{2}\right)  \sim2G\left(  1-r\right)  n
\end{array}
\]
for lazy random walks with reflection (in either strong or weak sense), which
surely must be known as well.

A \textbf{Motzkin prefix} of length $n$ is a ternary word of a total of $n$
$\nearrow$s, $\rightarrow$s and $\searrow$s in which no initial segment
contains more $\searrow$s than $\nearrow$s. \ (A dispersed Dyck prefix is a
Motzkin prefix with $\rightarrow$s prohibited except at zero height.) Counting
such words is summarized in \cite{S3-rw}. \ The associated (unweighted)
asymptotics remain open for arbitrary $p\neq q$ and $0\leq r\leq1$.

In particular, for $r=1/3$, we have%
\[%
\begin{array}
[c]{ccc}%
\mathbb{E}\left(  M_{n}\right)  /\sqrt{n}\sim\sqrt{\pi/3}\approx1.023, &  &
\mathbb{E}\left(  M_{n}^{2}\right)  /n\sim(4/3)G\approx1.221.
\end{array}
\]
Contrast this with%
\[%
\begin{array}
[c]{ccc}%
S_{0}=0, &  & S_{j}=\max\left\{  S_{j-1}+Y_{j},0\right\}
\end{array}
\]
in which $Y_{j}=-1$ when $j>0$ is divisible by $3$ and $\mathbb{P}\left\{
Y_{j}=1\right\}  =\mathbb{P}\left\{  Y_{j}=0\right\}  =1/2$ otherwise. \ Under
both circumstances, increments $\in\{1,-1,0\}$ are equiprobable, but the
timing of $-1$s for the latter is deterministic (similar to a traffic light
governing random traffic). \ We have \cite{F4-rw}%
\[%
\begin{array}
[c]{ccc}%
\mathbb{E}\left(  M_{n}\right)  /\sqrt{n}\sim\sqrt{\pi/12}\approx0.512, &  &
\mathbb{E}\left(  M_{n}^{2}\right)  /n\sim(1/3)G\approx0.305
\end{array}
\]
as $n=3k\rightarrow\infty$ via calculations resembling those in Sections 7 \&\ 8.

A\ different generalization allows the sequence $X_{1}$, $X_{2}$, \ldots,
$X_{n}$ to be dependent; more precisely, to follow a two-state Markov chain
with one-step transition matrix%
\[
\left(
\begin{array}
[c]{ll}%
\mathbb{P}\left\{  \left.  X_{i+1}=+1\right\vert \,X_{i}=+1\right\}  &
\mathbb{P}\left\{  \left.  X_{i+1}=-1\right\vert \,X_{i}=+1\right\} \\
\mathbb{P}\left\{  \left.  X_{i+1}=+1\right\vert \,X_{i}=-1\right\}  &
\mathbb{P}\left\{  \left.  X_{i+1}=-1\right\vert \,X_{i}=-1\right\}
\end{array}
\right)  =\left(
\begin{array}
[c]{cc}%
\alpha & \beta\\
\beta & \alpha
\end{array}
\right)
\]
and $\alpha+\beta=1$. \ Since $\alpha$ is usually understood to be $>1/2$, the
current step of the walk $S_{0}$, $S_{1}$, $S_{2}$, \ldots, $S_{n}$ tends to
move in the same direction as the preceding step (rather than changing
direction) and is called a \textbf{persistent random walk} \cite{Go-rw, Se-rw,
Ga-rw}. \ In the event $p=q$, we have \cite{Bh-rw}%
\[%
\begin{array}
[c]{ccc}%
\mathbb{E}\left(  M_{n}^{+}\right)  \sim\sqrt{\dfrac{\alpha}{\beta}}\left(
\sqrt{\dfrac{2n}{\pi}}-\dfrac{1}{2}\sqrt{\dfrac{\alpha}{\beta}}\right)  , &  &
\mathbb{V}\left(  M_{n}^{+}\right)  \sim\dfrac{\alpha}{\beta}\left(
1-\dfrac{2}{\pi}\right)  n.
\end{array}
\]
No one has investigated the cross-correlation $\mathbb{E}\left(  M_{n}%
^{+}\,M_{n}^{-}\right)  $, as far as is known, nor outcomes for the event
$p<q$.

Although their titles are similar to ours, the papers \cite{Y1-rw, Y2-rw} are
concerned not with moments of $M_{n}$ but rather moments of $S_{n}$.
\ Surveying $M_{n}$ results for dimensions $\geq2$, for both symmetric and
asymmetric cases, is left to some other researcher.

\section{Acknowledgements}

I am thankful to Irene Marzuoli, Richard Cowan, Niels Hansen and Kenneth S.
Williams for helpful correspondence.

\end{document}